\documentclass[12pt,psfig]{article}
\usepackage{graphicx,epsfig}
\newcommand{\si}{\sigma}
\newcommand{\ld}{\lambda}
\newcommand{\ba}{\begin{array}}
\newcommand{\ea}{\end{array}}
\newcommand{\ban}{\begin{eqnarray*}}
\newcommand{\ean}{\end{eqnarray*}}

\newcommand{\ol}{\overline}

\begin{document}

\baselineskip=17pt

\begin{center}

\vspace{-0.6 in} {\large \bf Trace class perturbation of closed linear relations $^{\dag}$}

 \vspace{0.3in}

Yuming Shi$^{a,{\ddag}}$ and  Yan Liu$^b$ \\

$^a$Department of Mathematics, Shandong University\\

Jinan, Shandong 250100, P. R. China\\

$^b$Department of Mathematics and Physics, Hohai University\\

Changzhou Campus 213022, P. R. China \\

\footnote{$^\dag$This research was supported by the NNSF of China (Grants 11571202).\\
\indent \ \ $^\ddag$The corresponding author.\\
\indent \indent Email addresses: ymshi@sdu.edu.cn(Y. Shi), yanliumaths@126.com (Y. Liu).}

\end{center}

{\bf Abstract.} This paper studies trace class perturbation of closed linear relations in Hilbert spaces.
The concept of  trace class perturbation of closed relations is introduced by orthogonal projections.
 Equivalent characterizations of  compact and trace class block operator matrices
are first given in terms of their elements, separately.  By using them, several equivalent and sufficient
 characterizations of trace class perturbation of closed linear relations
  are obtained.

\medskip

\medskip

\noindent{\bf 2010 AMS Classification}: 47A06, 47A55, 47B25, 47A10.

\bigskip

\noindent{\bf 1. Introduction}\medskip

In the classical operator theory, all the operators discussed are single-valued
(e.g., [8, 11, 20]). In the case that
an operator is not densely defined, its adjoint is multi-valued.
So it is always required that the operators are densely defined when one considers
 their adjoints in the classical operator theory.
In 1950, von Neumann introduced linear relations in order to study adjoints of non-densely defined
linear differential operators [10]. Since then, more and more multi-valued operators have been
found and then they have attracted a lot of  attention from mathematicians.
In 2003, Lesch and Malamud  studied symmetric linear differential expressions
whose  minimal  operators are  non-densely defined,
and whose maximal operators are multi-valued when the  differential expressions do not
satisfy the definiteness condition [9].
Recently, we found that minimal and maximal operators generated by symmetric
linear difference expressions are multi-valued or non-densely defined
in general even though the corresponding definiteness condition is satisfied [12, 15].
Obviously, the classical operator theory
is not available in this case. So it is very urgent for us to establish
the theory of multi-valued linear operators.

Multi-valued linear operators are often called linear relations (briefly, relations) or subspaces
of the related product spaces [1, 3, 10]. Linear relations include both
single-valued and multi-valued operators. Throughout the present paper,
an operator always means that it is single-valued for convenience.

Perturbation problems are very important  in pure
and applied mathematics. The classical perturbation theory of operators
has been studied for a long time, and some elegant results have been obtained (see [8, 11, 20]).
There have been some important progresses about perturbations of linear relations
made in the last decades. In 1998, Cross introduced a concept of  relatively compact perturbation
of  linear relations, and studied its some properties [4]. In 2009,
Azizov with his coauthors introduced concepts of
compact and finite rank perturbations of closed relations
in Hilbert spaces by orthogonal projections  (see Definition 2.3),
and gave some equivalent characterizations [2].
In 2014, Wilcox showed that five kinds of essential spectra of linear relations are
stable under relatively compact perturbation with some additional conditions and under
compact perturbation, separately  [21].
Motivated by the above works and the related existing results for linear operators,
the first author of the present paper
studied the stability of essential spectra of self-adjoint relations under compact
perturbation in 2016 [13]. She first studied the relationships among the operator
parts of the unperturbed relation, perturbed term and perturbed relation
using the decomposition of closed relations given by Arens [1]. Using these relationships,
she gave out invariance of self-adjointness and stability of essential spectra
of self-adjoint relations under compact perturbation [13].

In the present paper, we shall  study trace class perturbation of closed relations
in Hilbert spaces. Enlightened by the idea used
in the definitions of compact and finite rank perturbations of closed relations
in [2], we shall define  the trace class perturbation of closed relations by orthogonal
projections (see Definition 3.1). Then we shall study its characterizations, and
give out its several equivalent and sufficient
characterizations based on the research works in [2, 13, 14, 16, 18, 19].

The rest of the present paper is organized as follows.
In Section 2, some notations, basic concepts
and fundamental results about linear relations are introduced.
In particular, equivalent characterizations of compact and trace
class block operator matrices in terms of their elements are given, separately.
In Section 3, the concept of trace class perturbation of closed relations in Hilbert spaces is
introduced, and its several equivalent and sufficient characterizations are obtained.\medskip

\bigskip

\noindent{\bf 2. Preliminaries}
\medskip

In this section, we shall first list some notations and basic concepts, and
recall some fundamental results about linear relations, including
 resolvent set and spectrum of linear relations, and relationships
between them of  closed relations and those of their corresponding operator parts.
Then we shall recall the  concepts of finite rank and trace class operators and
their some properties. In addition, we shall give
equivalent characterizations of compact and trace
class block operator matrices in terms of their elements, separately,
which will be used in Section 3.
Finally, we shall introduce
the concepts of finite rank and compact perturbations of closed relations.

This section is divided into three subsections.

\medskip

\noindent{\bf 2.1. Some notations and basic concepts about linear relations}\medskip

In this subsection, we shall introduce  some notations and basic concepts of linear relations,
including closed, adjoint, Hermitian, and self-adjoint relations.

 By ${\mathbf R}$ and ${\mathbf C}$ denote the sets of the real and
complex numbers, respectively, throughout this paper.

Let $X, Y$, and $Z$ be linear spaces over a number field ${\mathbf K}$.
If $X$ is a normed space with norm $\|\cdot\|_X$ or an inner product space with inner
product $\langle \cdot, \cdot\rangle_X$,
the subscript $X$ will be omitted without confusion.
Denote ${\bar B}_X:=\{x\in X:\;\|x\|\le 1\}$ and
$\partial{\bar B}_X:=\{x\in X:\;\|x\|= 1\}$ if $X$ is a normed space.
If $X$ is an inner product space and $E\subset X$, by $E^\perp$
denote the orthogonal complement of $E$.

In the case that $X$ and $Y$ are topological linear spaces,
the topology of the product space $X\times Y$  is naturally induced by $X$ and $Y$.
Further, if $X$ and $Y$ are normed,
then the norm of $X\times Y$  is defined by
$$\|(x,y)\|=\left(\|x\|^2+\|y\|^2\right)^{1/2},\;(x,y)\in X\times Y.$$
Similarly, if $X$ and $Y$ are inner product spaces,  then the inner product
of $X\times Y$ is defined by
$$\langle(x_1,y_1),(x_2,y_2)\rangle=\langle x_1, x_2\rangle
+\langle y_1,y_2\rangle,\;\;(x_1,y_1),\;(x_2,y_2)\in X\times Y.$$

Every linear subspace $T\subset X\times Y$ is called a linear relation (briefly, relation or subspace) of $X
\times Y$. By $LR(X,Y)$ denote the set of all the linear relations of $X\times Y$.
In the case that $X=Y$, by $LR(X)$ denote $LR(X,Y)$ briefly.

Let $T\in LR(X,Y)$. The domain $D(T)$ and range $R(T)$ of $T$ are respectively defined by
$$\begin{array}{rrll}
D(T):&=&\{x\in X:\, (x,y)\in T \;{\rm for\; some}\;y\in Y\},\\[0.5ex]
R(T):&=&\{y\in Y:\, (x,y)\in T \;{\rm for\; some}\;x\in X\}.
\end{array}$$
Further, denote
$$T(x):=\{y\in Y:\,(x,y)\in T\},\;\;T^{-1}:=\{(y,x):\,(x,y)\in T\}.$$
It is evident that $T(0) = \{0\}$ if and only if
$T$ uniquely determines a linear operator
from $D(T)$ into $Y$ whose graph is $T$. For convenience,
a linear operator from $X$ to $Y$ will always be
identified with a subspace of $X\times Y$ via its graph.

In the case that  $X$ and $Y$ are topological linear spaces,  $T\in LR(X,Y)$ is said to be a closed relation
if ${\ol T}=T$, where ${\ol T}$ is  the closure of $T$.
By $CR(X,Y)$ denote the set of all the closed relations of  $X\times Y$.
By $CR(X)$ denote $CR(X,X)$ briefly.
It is evident that $T\in CR(X,Y)$ if and only if $T^{-1}\in CR(Y,X)$.

Let $S, T\in LR(X,Y)$ and $\alpha\in {\mathbf K}$. Define
$$\begin{array}{ccc}\alpha\, T :
= \{(x, \alpha\, y) : (x, y) \in T\},\\[0.5ex]
T + S := \{(x, y + z) : (x, y) \in T,\; (x, z)\in S\}.
\end{array}$$
If $T \cap S$ = \{(0, 0)\}, then denote
$$T\dot{+} S:=\{(x+u,y+v): \;(x,y)\in T,\; (u,v)\in S\}.$$
Further, in the case that $X$ and $Y$ are inner product spaces,
if $T$ and $S$ are orthogonal; that is,
$\langle(x, y),(u,v)\rangle = 0 $ for all
$(x, y) \in T$ and $(u, v) \in S$, then denote
$$T\oplus S:=T\dot{+} S.$$
Let $T\in LR(X,Y)$ and $S\in LR(Y,Z)$. The product of $T$ and $S$
is defined by
$$ST=\{(x,z)\in X\times Z:\;{\rm there\; exists}\; y \in Y\;
{\rm such\; that}\;(x,y)\in T\;{\rm and}\; (y,z)\in S\}.                             $$
Note that if $S$ and $T$ are operators, then $ST$ is also an operator.

Let $X$ be a Hilbert space. The adjoint of $T\in LR(X)$  is defined by
$$T^{*}:=\{(f,g)\in X^{2}:\;\langle g,x\rangle=\langle f,y\rangle\;{\rm for\;all}\;(x,y)\in T\}.$$
$T$ is said to be Hermitian in $X^{2}$ if $T\subset T^{*}$, and said to be self-adjoint
 in $X^{2}$ if $T=T^{*}$.
 \medskip

 \noindent{\bf Lemma 2.1  {\rm [16, Proposition 2.1]}.} Let $X$ and $Y$ be linear spaces, and $S, T\in LR(X,Y)$.
Then $S=(S-T)+T$  if and only if $D(S)\subset D(T)$ and $T(0)\subset S(0)$.\medskip

\noindent{\bf Lemma 2.2 {\rm [19, Proposition 3.1].}}
 Let $X$ and $Y$ be linear spaces, and $S,\,T\in LR(X,Y)$.
 If $S(0)\subset T(0)$ and $D(S)\subset D(T)$, then
$$T^{-1}-S^{-1}=T^{-1}(S-T)S^{-1}.     $$

Now, we shall  recall the definitions  of resolvent set and spectrum
of linear relations in complex Hilbert spaces.\medskip

 \noindent{\bf Definition 2.1  {\rm [6, 7, 14]}.} Let $X$ be a complex Hilbert space and $T\in LR(X)$.
The set $\rho(T):=\{\ld\in {\mathbf C}:\,(\ld I-T)^{-1}\;{\rm is\;  a\; bounded\;
 linear\; operator\; defined\; on}\;X\}$  is called the resolvent set of $T$, and
$\si(T):={\mathbf C}\setminus \rho(T)$ is called the spectrum of $T$.
\medskip

Let  $X$ be  a Hilbert space and $T\in CR(X)$.  Arens introduced  the following important decomposition [1]:
$$T=T_s\oplus T_\infty,                                                                                                                             \eqno (2.1) $$
where
$$T_\infty :=\{(0,g)\in X^2: (0,g)\in T\},\;\;
T_s :=T\ominus T_{\infty}.                                                                                                                     \eqno (2.2)$$
Then $T_s\in CR(X)$ is a linear operator, and $T_\infty\in CR(X)$.
So $T_s$ and $T_\infty$ are often called the operator
and pure multi-valued parts of $T$, respectively.
In addition, they satisfy the following properties [1]:
$$D(T_s)=D(T),\;R(T_s)\subset T(0)^\perp,\;
\;T_\infty =\{0\}\times T(0),                                                                                                                         \eqno (2.3)$$
and $D(T_s)$ is dense in $T^*(0)^\perp$.

We shall remark that this decomposition establishes an important bridge between
closed relations and operators. One can apply the properties of  the operator $T_s$
to study related problems about the closed relation $T$ in some cases (e.g., [13, 14, 16, 19]).
\medskip

\noindent{\bf Lemma 2.3 {\rm [13, Proposition 2.1]}.} Let $X$ be a Hilbert space
and $T\in CR(X)$ be Hermitian. Then
 $T$ is a self-adjoint relation in $X^2$ if and only if  $T_s$ is a self-adjoint operator  in $T(0)^\perp$.
\medskip

The necessity of the above result was given in [5, Page 26]. Throughout the present paper,
the resolvent set and spectrum of $T_s$ and $T_\infty$ mean those of
$T_s$ and $T_\infty$ restricted to $(T(0)^\perp)^2$ and $T(0)^2$, respectively.
\medskip

\noindent{\bf Lemma 2.4 {\rm [14, Proposition 2.1 and Theorem 2.1]}.}
Let $X$ be a complex Hilbert space, and $T\in CR(X)$ be Hermitian in $X^2$. Then
$$T_s=T\cap (T(0)^\perp)^2,\;\;T_\infty=T\cap T(0)^2,                                                        \eqno (2.4)$$
$T_s$ is a closed Hermitian operator in $T(0)^\perp$, $T_\infty$ is a closed
Hermitian relation in $T(0)^2$, and
$$\rho(T)=\rho(T_s),\;\; \si(T)=\si(T_s),\;\; \si(T_\infty)=\emptyset.                                              \eqno (2.5)   $$

\medskip

\noindent{\bf 2.2. Concepts of trace class and finite rank operators and their some properties}
\medskip

In this subsection, we recall the definitions of trace class and finite rank operators and give out its some properties.
For more discussions about it we refer to [20, Chaps. 6 and 7].
In particular, we shall give equivalent characterizations of compact and trace class block operator matrices in
terms of its elements, separately, which will be used in Section 3.

In this subsection, all the spaces discussed are Hilbert spaces.
Let $X$ and $Y$ be Hilbert spaces. For convenience, we shall introduce the following notations:
by $B(X,Y)$ denote all the bounded  operators from $X$ into $Y$,
by $EB(X,Y)$ denote all the bounded operators on $X$ into $Y$ (i.e., their domains are equal to the whole
space $X$),
by $D(X,Y)$ denote all the densely defined operators from $X$ into $Y$,
and denote $DB(X,Y):=D(X,Y)\cap B(X,Y)$.
Briefly, by  $B(X)$, $EB(X)$, $D(X)$, and $DB(X)$
denote  $D(X,X)$, $B(X,X)$, $EB(X,X)$, and $DB(X,X)$, respectively.

Let  $T$ be a compact operator on $X$ into $Y$. Then $T^*T$ is compact,
self-adjoint, and non-negative. Define $|T|=(T^*T)^{1/2}$ (see [20, Page 169]).
Then $|T|$ is compact, self-adjoint, and non-negative.
The non-zero eigenvalues of $|T|$ are called the singular values of $T$.
Let $\{s_n(T)\} $ denote the (possibly finite) non-increasing sequence of
the singular values of $T$ (every value counted
according its multiplicity as an eigenvalue of $|T|$). If
$$\sum_n s_n(T)<\infty,$$
then $T$ is said to be a trace class operator on $X$ into $Y$.
By $EB_1(X,Y)$ denote all the trace class operators on $X$ into $Y$.
By $EB_1(X)$ denote $EB_1(X,X)$ briefly, and
$T\in EB_1(X)$  is briefly called a trace class operator on $X$.

If $T\in EB(X, Y)$ satisfies $\dim R(T)<\infty$, then
$T$ is said to be a finite rank operator on $X$ into $Y$. It is evident that
if $T$ is a finite rank operator on $X$ into $Y$, then $T\in EB_1(X,Y)$.

The following result comes from [20, (a) and (c) of Theorem 7.8]. \medskip

\noindent{\bf Lemma 2.6.} Let $X$ and $Y$ be Hilbert spaces.
\begin{itemize}
\item[{\rm (i)}] If $S,T\in EB_1(X,Y) $, then $S+T\in EB_1(X,Y)$.

\item[{\rm (ii)}] If $T\in EB_1(X,Y)$ and $S\in EB(Y,Z)$, then $ST\in EB_1(X,Z)$. The corresponding assertion holds
for $T\in EB(X,Y)$ and $S\in EB_1(Y,Z)$.
\end{itemize}

\medskip

\noindent {\bf Lemma 2.7 {\rm [20, (a) of Theorem 4.14]}.}
  $T\in DB(X,Y)$ if and only if $T^*\in EB(Y,X)$.\medskip

Now, we study  characterizations of compact and trace class block operator matrices
and  give their equivalent characterizations in terms of their elements,
which will be used in the next section. We refer to [17] for more discussions about block operator matrices.
\medskip

\noindent{\bf Proposition 2.1.} Let $X$ and $Y$ be Hilbert spaces,  and $Q$
be an operator on $X\times Y$  into $X\times Y$ and  can be written as
$$Q
=\left(\begin{array}{cc} Q_{11}& Q_{12}\\Q_{21}&Q_{22}\end{array}\right) ,               \eqno (2.6)   $$
where $Q_{11}, Q_{21}, Q_{12}$ and $Q_{22}$ are operators on $X$ into $X$,
on $X$ into $Y$, on $Y$ into $X$, and on $Y$ into $Y$, respectively. Then
\begin{itemize}
\item[{\rm (i)}] $Q$ is compact on $X\times Y$ if and only if $Q_{ij}$, $1\le i,j\le 2$, are compact
on their corresponding spaces, respectively;
\item[{\rm (ii)}] $Q\in EB_1(X\times Y)$  if and only if
$Q_{ij}$, $1\le i,j\le 2$, are all trace class operators
on their corresponding spaces, respectively.
\end{itemize}

\noindent{\bf Proof.} The assertion (i) can be easily verified by the definition of compact operators (see [20, Page 130])
 and (2.6).

Now, we show that the assertion (ii) holds.
We shall first show that its necessity holds. Suppose that  $Q\in EB_1(X\times Y)$.
Then
$$\sum s_n(Q)<\infty.                                                                                                                                    \eqno (2.7)   $$
It follows from (2.6) that
$$Q^*Q=\left(\begin{array}{cc}Q_{11}^* Q_{11}+Q_{21}^*Q_{21}& Q_{11}^*Q_{12}+Q_{21}^*Q_{22}
\\Q_{12}^*Q_{11}+Q_{22}^*Q_{21}&Q_{12}^*Q_{12}+Q_{22}^*Q_{22}\end{array}\right).                              \eqno (2.8)$$
Denote
$$Q_i:=Q_{1i}^* Q_{1i}+Q_{2i}^*Q_{2i},\;\;i=1, 2.                                                                                \eqno (2.9)  $$
It can be easily verified that $Q_1$ and $Q_2$ are compact, self-adjoint and non-negative operators
on $X$ and on $Y$, respectively. Denote
$$P_i:=Q_i^{1/2},\;\;i=1, 2.                                                                                                                      \eqno (2.10)  $$
We shall only show that $Q_{11}\in EB_1(X)$. With a similar argument,
one can show that the others hold.  Let $\ld$ be any non-zero eigenvalue of $P_1$. Then $\ld^2$ is an eigenvalue of $Q_1$,
and so is an eigenvalue of $Q^*Q$ by (2.8). Hence, $\ld$ is a singular value of $Q$. By (2.7) we get
that $P_1\in EB_1(X)$. In addition, it follows from (2.9)  and (2.10) that
$Q_{11}^*Q_{11}\le Q_1=P_1^2$. Then
$$\langle(P_1^2-Q_{11}^*Q_{11})(x),x\rangle\ge 0, \;\;\forall\,x\in X,$$
which implies that
$$\|Q_{11}(x)\|\le \|P_1(x)\|,\;\; \forall\,x\in X.                                                                                       \eqno (2.11)    $$
Further, by [20, Theorem 7.7] and (2.11) we have that
$$s_1(Q_{11})=\|Q_{11}\|\le \|P_1\|=s_1(P_1),                                                                         $$
$$\begin{array}{rrll}
&&s_{n+1}(Q_{11})=\inf\limits_{\stackrel{x_k\in X}{1\le k\le n}}\sup
\{\|Q_{11}(x)\|:\,x\bot x_1,\ldots,x_n \;
{\rm for\;all}\;x\in \partial {\bar B}_{X}\}\\
&\le & \inf\limits_{\stackrel{x_k\in X}{1\le k\le n}}\sup
\{\|P_1(x)\|:\,x\bot x_1,\ldots,x_n \;
{\rm for\;all}\;x\in \partial {\bar B}_{X}\}\\
&=& s_{n+1}(P_1),\;\;n\ge 1.
\end{array}                                                                                                                               $$
This yields that $\sum s_n(Q_{11})\le \sum s_n(P_1) <\infty$. Therefore, $Q_{11}\in EB_1(X)$,
and consequently the necessity of the assertion (ii) holds.

Next, we shall consider the sufficiency of the assertion (ii). Suppose that
$Q_{ij}$, $1\le i,j\le 2$, are all trace class operators
on their corresponding spaces, respectively.
Then
 $$\sum s_n(Q_{ij}) <\infty,\;\;1\le i,j\le 2.                                                                                  \eqno (2.12)  $$
Again by [20, Theorem 7.7] we have that
$$\begin{array}{rrll}
&&s_1(Q)=\|Q\| =\sup\{\|Q(x,y)\|:\, (x,y)\in \partial {\bar B}_{X\times Y}\}\\ [0.5ex]
&=& \sup\{\|Q(x,y)\|:\,(x,y)\in  {\bar B}_{X\times Y}\}\\ [0.5ex]
&=&\sup\{\Big(\|Q_{11}(x)+Q_{12}(y)\|^2+\|Q_{21}(x)+Q_{22}(y)\|^2\Big)^{\frac{1}{2}}:\,(x,y)\in  {\bar B}_{X\times Y}\}\\[0.5ex]
&\le& \sup\{\|Q_{11}(x)\|+\|Q_{12}(y)\|+\|Q_{21}(x)\|+\|Q_{22}(y)\|:\,(x,y)\in  {\bar B}_{X\times Y}\}\\[0.5ex]
&\le &  \sum_{i=1}^2\sup\{\|Q_{i1}(x)\|:\,x\in {\bar B}_X\}+ \sum_{i=1}^2\sup\{\|Q_{i2}(y)\|:\,y\in {\bar B}_Y\}\\[0.5ex]
&=& \sum_{i,j=1}^2\|Q_{ij}\|= \sum_{i,j=1}^2s_1(Q_{ij}).
 \end{array}                                                                                                                         \eqno (2.13)  $$
Similarly, for any given $n\ge1$ we get that
$$\begin{array}{rrll}
&&s_{n+1}(Q)=\inf\limits_{\stackrel{f_k\in X\times Y}{1\le k\le n}}\sup
\{\|Q(x,y)\|:\,(x,y)\bot f_1, \ldots, f_n \;
{\rm for\;all}\;(x,y)\in \partial {\bar B}_{X\times Y}\}\\
&=&\inf\limits_{\stackrel{f_k\in X\times Y}{1\le k\le n}}\sup
\{\|Q(x,y)\|:\,(x,y)\bot f_1, \ldots, f_n \; {\rm for\;all}\;(x,y)\in {\bar B}_{X\times Y}\}\\
&\le & \inf\limits_{\stackrel{f_k\in X\times Y}{1\le k\le n}}
 \sup\{\sum_{i=1}^2\|Q_{i1}(x)\|+\sum_{i=1}^2\|Q_{i2}(y)\|:\,
 (x,y)\bot f_1, \ldots, f_n \; {\rm for\;all}\;(x,y)\in {\bar B}_{X\times Y}\}\\
&\le & \gamma_1+\gamma_2+\gamma_3+\gamma_4,
\end{array}                                                                                                                              \eqno (2.14)  $$
where
$$\gamma_i= \inf\limits_{\stackrel{f_k\in X\times Y}{1\le k\le n}}
 \sup\{\|Q_{i1}(x)\|:\,(x,y)\bot f_1, \ldots, f_n \; {\rm for\;all}\;(x,y)\in {\bar B}_{X\times Y}\},
 \;\;i=1,2$$
 $$\gamma_{i+1}=\inf\limits_{\stackrel{f_k\in X\times Y}{1\le k\le n}}
 \sup\{\|Q_{i2}(y)\|:\,(x,y)\bot f_1, \ldots, f_n \; {\rm for\;all}\;(x,y)\in {\bar B}_{X\times Y}\},
 \;\;i=1,2.$$
For any given $x_1,\ldots,x_n \in X$,   we obtain that
$$\begin{array}{rrll}
&&
\gamma_1\le \sup\{\|Q_{11}(x)\|:\,(x,y)\bot (x_1,0), \ldots, (x_n,0) \; {\rm for\;all}\;(x,y)\in {\bar B}_{X\times Y}\}\\[0.6ex]
&= & \sup\{\|Q_{11}(x)\|:\,x\bot x_1, \ldots, x_n \; {\rm for\;all}\;x\in {\bar B}_{X}\},
\end{array}  $$
which, together with  [20, Theorem 7.7], implies that
$$\gamma_1\le s_n(Q_{11}).                                                                                                        \eqno (2.15) $$
With a similar argument to the above, one can get that
$$\gamma_2\le s_n(Q_{21}), \;\gamma_3\le s_n(Q_{12}),\; \gamma_4\le s_n(Q_{22}).                      \eqno (2.16)$$
It follows from (2.14), (2.15) and (2.16) that
$$s_{n+1}(Q)\le  \sum_{i,j=1}^2s_{n+1}(Q_{ij}),\;\;n\ge 1.                                                                                                                        $$
This, together with (2.12)  and  (2.13), yields that (2.7) holds.
Therefore, $Q\in EB_1(X\times Y)$, and consequently
the sufficiency of the assertion (ii) holds.

The whole proof is complete.
\bigskip

\noindent{\bf 2.3. Concepts of finite rank and compact perturbations of closed relations}
\medskip

In this subsection, we shall recall the definitions of of finite rank and compact perturbations
of closed relations. We refer to [2] for more discussions.\medskip

Let $X$ be a Hilbert space, and $M$  be a closed subspace of $X$.
By $P_M^X$ denote the orthogonal projection from $X$ onto $M$.
The superscript $X$ of $P_M^X$ is omitted without confusion.\medskip

\noindent{\bf Definition 2.3.} Let $X$ and $Y$ be  Hilbert spaces,  $S, T\in CR(X,Y)$, and
$$P_T:\,X\times Y\to T,\;\;P_S:\,X\times Y\to S                                                    $$
 be orthogonal projections.
\begin{itemize}
\item[{\rm (1)}] $T$ is said to be a finite rank perturbation
of $S$ in $X\times Y$ if $P_T-P_S$ is a finite rank operator on  $X\times Y$.

\item[{\rm (2)}] $T$ is said to be a compact perturbation of $S$ in $X\times Y$
if $P_T-P_S$ is a compact operator on  $X\times Y$.
\end{itemize}
\medskip

It is evident that if $T$ is a finite rank perturbation
of $S$ in $X\times Y$, then $T$ is a compact perturbation of $S$ in $X\times Y$.

\bigskip

\noindent{\bf 3. Concept of trace class perturbation of closed relations and its characterizations}
\medskip

In this section, we shall pay our attention to  trace class perturbation of closed relations and its characterizations.
We shall first introduce the definition of trace class perturbation of closed relations, and then give out
its several equivalent and sufficient characterizations.

Throughout this section, we always assume that  $X$ and $Y$ are complex Hilbert spaces.

Enlightened by the definitions of compact and finite rank perturbations of closed relations
given in [2] (see Definition 2.3), we introduce the following definition of
trace class perturbation of closed relations:
\medskip

\noindent{\bf Definition 3.1.} Let  $S, T\in CR(X,Y)$, and
$P_S$ and $P_T$ be orthogonal projections from $X\times Y$ onto $S$ and $T$, respectively.
Then $T$ is said to be a trace class perturbation of $S$ in $X\times Y$
if $P_T-P_S\in EB_1(X\times Y)$.
\medskip

It is evident that if $T$ is  a trace class perturbation of $S$ in $X\times Y$, then
$T$ is  a compact  perturbation of $S$ in $X\times Y$; and
if $T$ is  a finite rank perturbation of $S$ in $X\times Y$, then
$T$ is  a trace class  perturbation of $S$ in $X\times Y$ by
Definitions 2.3 and 3.1.\medskip

\noindent {\bf Lemma 3.1.} Let $S, T\in CR(X,Y)$ and $A\in B(X,Y)$ with $D(S)\cup D(T)\subset D(A)$.
 Then the following inequalities hold:
$$\frac{1}{\gamma}\|P_{T-A}-P_{S-A}\| \le \|P_{T}-P_{S}\| \le \gamma \|P_{T-A}-P_{S-A}\|,         $$
where $\gamma:=2(1+\|A\|^2)$, and $P_T$, $P_S$, $P_{T-A}$, and $P_{S-A}$ are
orthogonal projections from $X\times Y$ onto $T$,  $S$, $T-A$, and $S-A$, respectively.
\medskip

Motivated by [2, Lemma 4.1], we give out the above result. Note that it is required that $A\in EB(X,Y)$
in [2, Lemma 4.1].  However, this assumption can be weakened
 by $A\in B(X,Y)$ with $D(S)\cup D(T)\subset D(A)$. The proof of
 Lemma 3.1 is similar to that of [2, Lemma 4.1], and
so its details are omitted.\medskip

\noindent{\bf Proposition 3.1.} Let $S, T\in CR(X,Y)$ and $A\in DB(X,Y)$ with $D(S)\cup D(T)\subset D(A)$.
Then $T$ is a trace class perturbation of $S$ in $X\times Y$
if and only if $T-A$ is a trace class perturbation of $S-A$ in $X\times Y$.\medskip

\noindent{\bf Proof.} It can be easily verified that $S-A, T-A\in CR(X,Y)$
by the closedness of $S$ and $T$, and the boundedness of $A$.

By Definition 3.1 it suffices to show that $P_T-P_S\in EB_1(X\times Y)$
if and only if $P_{T-A}-P_{S-A}\in EB_1(X\times Y)$.
So it suffices to show that
the following inequalities hold:
$$\frac{1}{\gamma}s_n(P_{T-A}-P_{S-A})\leq s_n(P_{T}-P_{S})
\leq \gamma s_n(P_{T-A}-P_{S-A}),\;n\ge 1,                                                   \eqno (3.1)$$
where $\gamma:=2(1+\|A\|^2)$, and $\{s_n(P_{T}-P_{S})\}$
and $\{s_n(P_{T-A}-P_{S-A})\}$
are the non-increasing sequences of
the singular values of $P_{T}-P_{S}$ and $P_{T-A}-P_{S-A}$, respectively.
Since $A(0)=\{0\}\subset T(0)\cap S(0)$ and $D(S)\cup D(T)\subset D(A)$, it follows from Lemma 2.1 that
$$T-A+A=T,\;\;S-A+A=S.                                                $$
Hence, it is only needed for us to show that the first inequality in (3.1),
namely
$$s_n(P_{T-A}-P_{S-A})\leq \gamma s_n(P_{T}-P_{S}),\;\;n\ge 1,                                                                \eqno (3.2)       $$
holds because the second inequality in (3.1) follows by replacing
$T, S$ and $A$ with $T-A, S-A$ and $-A$, respectively.

By [20, Theorem 7.7 and Exercise 7.2] and Lemma 3.1,
we have that
$$\frac{1}{\gamma}s_1(P_{T-A}-P_{S-A})=\frac{1}{\gamma}\|P_{T-A}-P_{S-A}\|
\leq \|P_{T}-P_{S}\|=s_1(P_{T}-P_{S}),                                                                                             \eqno (3.3) $$
$$\begin{array}{rrll}
&&s_n(P_{T}-P_{S})\\
&=&\min\limits_{\stackrel{\varphi_j\in X\times Y}{1\le j\le n-1}}{\rm max}
\{\|(P_{T}-P_{S})(\psi)\|:\,\psi\bot \varphi_1,\ldots,\varphi_{n-1} \;
{\rm for\;all}\;\psi\in \partial {\bar B}_{X\times Y}\},
\end{array}                                                                                                                                \eqno (3.4)$$
$$\begin{array}{rrll}
&&s_n(P_{T-A}-P_{S-A})\\
&=&\min\limits_{\stackrel{\tilde{\varphi}_j\in X\times Y}{1\le j\le n-1}}
{\rm max} \{\|(P_{T-A}-P_{S-A})(\tilde{\psi})\|:\,\tilde{\psi}\bot \tilde{\varphi}_1,\ldots,\tilde{\varphi}_{n-1}
 \;{\rm for\;all}\;\tilde{\psi}\in \partial {\bar B}_{X\times Y}\}
 \end{array}                                                                                                                                    \eqno (3.5)$$
for all $n\ge 2$. Hence, (3.2) holds for $n=1$ by (3.3).

Now, we shall show that (3.2) holds for any given $n\ge 2$.
It follows from (3.4) that there exist
$\varphi_1=(x_1,y_1),\ldots,\varphi_n=(x_n,y_n)\in X\times Y$ such that
$$s_n(P_{T}-P_{S})={\rm max}\{\|(P_{T}-P_{S})\psi\|:\,
\psi\bot \varphi_1,\ldots,\varphi_{n-1} \;
{\rm for\;all}\;\psi\in \partial {\bar B}_{X\times Y}\}.                                                                           \eqno (3.6)     $$
By the assumption that $A\in DB(X,Y)$ and  by Lemma 2.7
we have that $A^*\in EB(Y, X)$. Let
$$\begin{array}{rrll}
M&:=&L\{\varphi_1,\ldots,\varphi_{n-1}\},\\[1ex]
\tilde{M}&:=&L\{(x_1+A^*y_1,y_1),...,(x_{n-1}+A^*y_{n-1},y_{n-1})\}.
\end{array}$$
Then it follows from (3.6) and (3.5) that
$$\begin{array}{rrll}
s_n(P_{T}-P_{S})
&=&{\rm max}\{\|(P_{T}-P_{S})\psi\|:\,\psi\in \partial {\bar B}_{{M}^\perp}\}, \\
s_n(P_{T-A}-P_{S-A})&\leq& {\rm max}\{\|(P_{T-A}-P_{S-A})\tilde{\psi}\|:\,
\tilde{\psi}\in \partial {\bar B}_{{\tilde{M}}^\perp}\}.
\end{array}                                                                                                                                     \eqno (3.7) $$
With the help of
$$\|P_{T}-P_{S}\|={\rm max}\{\sup\limits_{\omega\in \partial {\bar B}_T}{\rm d}(\omega,S),
\sup\limits_{\eta\in \partial {\bar B}_S}{\rm d}(\eta,T)\}                                       $$
(cf. [8, Page 198]), one can easily get that
$$\max\limits_{\psi\in \partial {\bar B}_{M^\perp}}\|(P_{T}-P_{S})\psi\|=
{\rm max}\{\sup\limits_{\omega\in \partial {\bar B}_{T\cap M^\perp}}
{\rm d}(\omega,S\cap {M}^\perp),
\sup\limits_{\eta\in \partial {\bar B}_{S\cap M^\perp}}
{\rm d}(\eta,T\cap {M}^\perp)\},$$
$$\begin{array}{rrll}
&&\max\limits_{\tilde{\psi}\in \partial {\bar B}_{{\tilde{M}}^\perp}}
\|(P_{T-A}-P_{S-A})\tilde{\psi}\|\\
&=&
{\rm max}\{\sup\limits_{\tilde{\omega}\in \partial {\bar B}_{(T-A)\cap {\tilde{M}}^\perp}}
{\rm d}(\tilde{\omega},(S-A)\cap {\tilde{M}}^\perp),
\sup\limits_{\tilde{\eta}\in \partial {\bar B}_{(S-A)\cap {\tilde{M}}^\perp}}
{\rm d}(\tilde{\eta},(T-A)\cap {\tilde{M}}^\perp)\},
\end{array}  $$
which, together with (3.7), implies that
$$s_n(P_{T}-P_{S})=
{\rm max}\{\sup\limits_{\omega\in \partial {\bar B}_{T\cap {M}^\perp}}
{\rm d}(\omega,S\cap {M}^\perp),
\sup\limits_{\eta\in \partial {\bar B}_{S\cap {M}^\perp}}
{\rm d}(\eta,T\cap {M}^\perp)\},                                                                                      \eqno (3.8)$$
$$\begin{array}{rrll}
&&s_n(P_{T-A}-P_{S-A})\\
& \leq &
{\rm max}\{\sup\limits_{\tilde{\omega}\in \partial {\bar B}_{(T-A)\cap {\tilde{M}}^\perp}}
{\rm d}(\tilde{\omega},(S-A)\cap {\tilde{M}}^\perp),
\sup\limits_{\tilde{\eta}\in \partial {\bar B}_{(S-A)\cap {\tilde{M}}^\perp}}
{\rm d}(\tilde{\eta},(T-A)\cap {\tilde{M}}^\perp)\}.
\end{array}                                                                                                                           \eqno (3.9) $$

In order to show that (3.2) holds, we shall first show that
$$\sup\limits_{\tilde{\omega}\in \partial {\bar B}_{(T-A)\cap {\tilde{M}}^\perp}}
{\rm d}(\tilde{\omega},(S-A)\cap {\tilde{M}}^\perp)
\leq \gamma s_n(P_{T}-P_{S}).                                                                                          \eqno (3.10)$$
In the case that $(T-A)\cap {\tilde{M}}^\perp=\{(0,0)\}$, we have that
$$\sup\limits_{\tilde{\omega}\in \partial {\bar B}_{(T-A)\cap {\tilde{M}}^\perp}}
{\rm d}(\tilde{\omega},(S-A)\cap  {\tilde{M}}^\perp)=0.$$
So (3.10) holds obviously in this case.
In the other case that $(T-A)\cap {\tilde{M}}^\perp\neq\{(0,0)\}$,
for any $\tilde{\omega}\in \partial {\bar B}_{(T-A)\cap {\tilde{M}}^\perp}$,
there exists $(u,v)\in T$ such that $\tilde{\omega}=(u,v-Au)$ with
$$\|u\|^2+\|v-Au\|^2=1.                                                                                                          \eqno (3.11)$$
We claim that $(u,v)\in {M}^\perp. $
In fact, for any $(x,y)\in M$, it can be easily verified that
$(x+A^*y,y)\in \tilde{M}$.
Then we have that
$$\langle (u,v-Au),(x+A^*y,y)\rangle=0,$$
which yields that $\langle u,x \rangle+ \langle v,y \rangle=0,$
namely $(u,v)\bot (x,y)$, and so $(u,v)\in {M}^\perp$.
Then $(r^{-1}u,r^{-1}v)\in \partial {\bar B}_{T\cap {M}^\perp}$, where
$r:=\left(\|u\|^2+\|v\|^2\right)^{1/2}>0$.
Therefore, for any given $\delta>s_n(P_{T}-P_{S})$, by (3.8)
there exists $(u',v')\in S\cap {M}^\perp$ such that
$$\|r^{-1}u-u'\|^2+\|r^{-1}v-v'\|^2<\delta^2,$$
i.e.,
$$\|u-ru'\|^2+\|v-rv'\|^2<r^2\delta^2.                                                                                                       \eqno (3.12)$$
Set
$$\psi_0:=r(u',v'-Au').$$
It can be easily verified that $\psi_0\in (S-A)\cap {\tilde{M}}^\perp$.
By (3.12) we have that
$$\begin{array}{rrll}
\|\tilde{\omega}-\psi_0\|^2
&=&\|u-ru'\|^2+\|v-r v'-A (u-ru')\|^2\\[1ex]
&\leq &\|u-ru'\|^2+2\|v-r v'\|^2+2\|A\|^2 \|u-ru'\|^2\\[1ex]
&\leq &2(1+\|A\|^2)(\|u-ru'\|^2+\|v-r v'\|^2)<\gamma r^2\delta^2.
\end{array}$$
In addition, it follows from (3.11) that
$$r^2=\|u\|^2+\|v-Au+Au\|^2\leq \|u\|^2+2\|v-Au\|^2+2\|A\|^2 \|u\|^2\leq \gamma.$$
So we get that
$$\|\tilde{\omega}-\psi_0\|< \gamma \delta,$$
which implies that
$$d(\tilde{\omega},(S-A)\cap{\tilde{M}}^\perp)\le \|\tilde{\omega}-\psi_0\|< \gamma \delta. $$
Therefore, by the arbitrariness of $\delta$ we get that
$$d(\tilde{\omega},(S-A)\cap{\tilde{M}}^\perp)\le \gamma s_n(P_{T}-P_{S}),$$
which yields that (3.10) holds.

With a similar argument, one can show that
$$\sup\limits_{\tilde{\eta}\in \partial {\bar B}_{(S-A)\cap {\tilde{M}}^\perp}}
{\rm d}(\tilde{\eta},(T-A)\cap {\tilde{M}}^\perp)\leq \gamma s_n(P_{T}-P_{S}),$$
which, together with (3.9) and (3.10), yields that (3.2) holds.
This completes the proof.\medskip

The following two results are direct consequences of Proposition 3.1.\medskip

\noindent{\bf Theorem 3.1.} Let $S, T\in CR(X,Y)$.
Then $T$ is a trace class perturbation of $S$ in $X\times Y$
if and only if $T-A$ is a trace class perturbation of $S-A$ in $X\times Y$ for some (and hence for all)
$A\in DB(X,Y)$ with $D(S)\cup D(T)\subset D(A)$.\medskip

\noindent{\bf Theorem 3.2.} Let $S, T\in CR(X)$.
Then $T$ is a trace class perturbation of $S$ in $X^2$
if and only if $T-\lambda I$ is a trace class perturbation of $S-\lambda I$ in $X^2$
for some (and hence for all) $\lambda\in {\mathbf{C}}$.\medskip

For $S,T\in CR(X,Y)$, we introduce the following set:
$$\Gamma (S, T):=\{A\in DB(X,Y)\; {\rm with}\; D(S)\cup D(T)
\subset D(A):\, (S-A)^{-1}, (T-A)^{-1}\in DB(Y,X)\}.$$

\noindent {\bf Remark 3.1.} $\Gamma (S, T)$ may be empty in some cases. For example,
let $Y={\mathbf C}^2$, $X=\{(x_1,0):\;x_1\in {\mathbf C}\}\subset Y$,
and $S, T\in CR(X,Y)$ be defined by $S(x_1,0)=\{(0,0)\}$
and $T(x_1,0)=\{(x_1,0)\}$  for any $(x_1,0)\in X$, respectively.
For every $A\in DB(X,Y)$, $(S-A)^{-1}=-A^{-1}\notin DB(Y,X)$
since $\dim R(A)<2$. Hence, $\Gamma (S, T)=\emptyset$ in this case.\medskip

\noindent{\bf Proposition 3.2.} Let $S, T\in CR(X,Y)$
and $\Gamma (S, T)\ne \emptyset$. And
let $A\in \Gamma (S, T)$.
Then $T$ is a trace class perturbation of $S$  in $X\times Y$ if and only
if $(T-A)^{-1}-(S-A)^{-1}\in EB_1(Y, X)$.\medskip

\noindent{\bf Proof.} For convenience, set $W:=(T-A)^{-1}-(S-A)^{-1}$.

By Proposition 3.1,  $T$ is a trace class perturbation of $S$ in $X\times Y$ if and only if
$T-A$ is a trace class perturbation of $S-A$ in $X\times Y$, namely
$P_{T-A}-P_{S-A}\in EB_1(X\times Y)$  by Definition 3.1.

Observe that for every $(x,y)\in X\times Y$,
$$P_{T-A}(x,y)=(x',y')$$
if and only if
$$P_{(T-A)^{-1}}(y,x)=(y',x').$$
$P_{S-A}$ and $P_{(S-A)^{-1}}$ have the same relation as the above.
Therefore,  $T$ is a trace class perturbation of $S$ in $X\times Y$  if and only if
$P_{(T-A)^{-1}}-P_{(S-A)^{-1}}\in EB_1(Y\times X)$.

In addition, by [2, Corollary 2.2], $P_{(T-A)^{-1}}-P_{(S-A)^{-1}}$ can be decomposed as
$$P_{(T-A)^{-1}}-P_{(S-A)^{-1}}
=\left(\begin{array}{cc} P_{11}& P_{12}\\P_{21}&P_{22}\end{array}\right),                  \eqno (3.13)    $$
where
$$P_{11}=F(T)-F(S),\; \;  P_{21}=(T-A)^{-1}F(T)-(S-A)^{-1}F(S),                                  \eqno (3.14)$$
$$P_{12}=(T^*-A^*)^{-1}H(T)-(S^*-A^*)^{-1}H(S),                                             \eqno (3.15)$$
$$P_{22}=(T-A)^{-1}(T^*-A^*)^{-1}H(T)-(S-A)^{-1}(S^*-A^*)^{-1}H(S),       \eqno (3.16)$$
while
$$F(T)=\left[I+(T^*-A^*)^{-1}(T-A)^{-1}\right]^{-1},\;
H(T)=\left[I+(T-A)^{-1}(T^*-A^*)^{-1}\right]^{-1}.                                                    \eqno (3.17)$$
Therefore,  $T$ is a trace class perturbation of $S$ in $X\times Y$  if and only if
$P_{ij}$, $1\le i,j\le 2$, are all trace class operators
on their corresponding spaces, respectively, by Proposition 2.1.

{\bf Necessity.} Suppose that $T$ is a trace class perturbation of $S$ in $X\times Y$.
Then $P_{ij}$, $1\le i,j\le 2$, are all trace class operators
on their corresponding spaces, respectively.
It follows from (3.14) and (3.17) that
$$W
=-\left[(S-A)^{-1}P_{11}-P_{21}\right]F(T)^{-1}.                   \eqno (3.18)   $$
So $W\in EB_1(Y, X)$ by Lemma 2.6 and by the fact that $(S-A)^{-1}$ and
$F(T)^{-1}$ are bounded.

{\bf Sufficiency.} Suppose that $W\in EB_1(Y, X)$.
Then $W^*\in EB_1(X, Y)$ by [20, Theorem 7.6].
In order to show that $T$ is a trace class perturbation of $S$ in $X\times Y$, it suffices to show that
 $P_{ij}$, $1\le i,j\le 2$,
are trace class operators on their corresponding spaces, respectively, by the above discussion.

It follows from (3.14), (3.17), and (3.18)  that
$$P_{11}=F(T)L_1F(S),\;P_{21}=(S-A)^{-1}P_{11}+WF(T),                                    \eqno (3.19)  $$
where
$$\begin{array}{rrll}
L_1 &=&(S^*-A^*)^{-1}(S-A)^{-1}-(T^*-A^*)^{-1}(T-A)^{-1}\\
&=&-(S^*-A^*)^{-1}W-W^*(T-A)^{-1}.
\end{array}                                                                                                                       \eqno (3.20)$$
Note that $L_1\in EB_1(Y)$  again
by Lemma 2.6 and by the fact that $(S^*-A^*)^{-1}$ and $(T-A)^{-1}$ are bounded.
Hence, it follows from (3.19) that $P_{11}\in EB_1(Y)$ and $P_{21}\in EB_1(Y, X)$
by Lemma 2.6  and by the fact that $F(T)$, $F(S)$, and $(S-A)^{-1}$ are bounded.

With a similar argument, interchanging $(T-A)^{-1}$ and $(S-A)^{-1}$ with $(T^*-A^*)^{-1}$
 and $(S^*-A^*)^{-1}$, separately, in (3.13) and (3.20),
 one can get that $P_{12}$, $H(T)-H(S)$, and the operator
$$L_2=(S-A)^{-1}(S^*-A^*)^{-1}-(T-A)^{-1}(T^*-A^*)^{-1}$$
are all trace class operators on their corresponding spaces, respectively.
Hence, the operator $L_3=L_2H(T)\in EB_1(X)$.
By (3.15) and (3.16) it can be easily verified that
$$P_{22}=(S-A)^{-1}(S^*-A^*)^{-1}\left[H(T)-H(S)\right]-L_3.$$
Thus, $P_{22}\in EB_1(X)$. Therefore,
$T$ is a trace class perturbation of $S$ in $X\times Y$.

The entire proof is complete.

\medskip

By Proposition 3.2 one can get the following two results:\medskip

\noindent{\bf Theorem 3.3.} Let $S, T\in CR(X,Y)$ and $\Gamma (S, T)\ne \emptyset$.
Then $T$ is a trace class perturbation of $S$ in $X\times Y$
if and only if $(T-A)^{-1}-(S-A)^{-1}\in EB_1(Y, X)$
 for some (and hence for all) $A\in \Gamma (S, T)$.\medskip

\noindent{\bf Theorem 3.4.} Let $S, T\in CR(X)$ and
$\rho(S)\cap \rho(T)\ne \emptyset$.
Then $T$ is a trace class perturbation of $S$ in $X^2$
if and only if $(T-\lambda I)^{-1}-(S-\lambda I)^{-1}\in EB_1(X)$
 for some (and hence for all) $\lambda\in \rho(S)\cap \rho(T)$.\medskip

Next, we shall give out other several equivalent characterizations of trace class perturbation
in terms of  the operator parts of $T$ and $S$ under some additional conditions.
\medskip

\noindent{\bf Theorem 3.5.} Let $S, T\in CR(X)$ satisfy that
 $S(0)=T(0)$ and $D(S)\cup D(T)\subset T(0)^\perp$.
Then $T$ is a trace class perturbation of $S$ in $X^2$ if and only if
$T_s$ is a trace class perturbation of $S_s$ in $(S(0)^\perp)^2$,
where $T_s$ and $S_s$ are the operator parts of $T$ and $S$, respectively.
\medskip

\noindent{\bf Proof.} By the assumption that $S(0)=T(0)$ and
$D(S)\cup D(T)\subset S(0)^\perp$ and by (2.2) and (2.3)
 we have  that
$$T_s=T\cap (S(0)^\perp)^2,\;S_s=S\cap (S(0)^\perp)^2,\;
T_\infty=S_\infty=\{0\}\times S(0).                                                                                   \eqno (3.21)$$
Fix any $(x,y)\in X^2$. There exist $x_1,y_1
\in S(0)^\perp$ and $x_2,y_2\in S(0)$
such that $(x,y)=(x_1,y_1)+(x_2,y_2)$.
It follows from (3.21) that
$$P_T(x,y)=P_{T_s}(x_1,y_1)+(0,y_2),\;
P_S(x,y)=P_{S_s}(x_1,y_1)+(0,y_2),$$
which implies that
$$(P_T-P_S)(x,y)=(P_{T_s}-P_{S_s})(x_1,y_1).                 $$
Therefore, the result of the theorem holds by Definition 3.1.
This completes the proof.
\medskip

\noindent{\bf Corollary 3.1.} Let $S, T\in CR(X)$ be Hermitian with
 $S(0)=T(0)$. Then the result of Theorem 3.5 holds.
\medskip

\noindent{\bf Proof.} Since $S$ and $T$ are Hermitian with $S(0)=T(0)$,
we get that $D(S)\subset S(0)^\perp$ and $D(T)\subset T(0)^\perp=S(0)^\perp$
by (2.4). So all the assumptions of Theorem 3.5 are satisfied, and consequently
the result of Theorem 3.5  holds. The proof is complete.
\medskip

\noindent{\bf Corollary 3.2.} Let $S, T\in LR(X)$ be self-adjoint
with  $D(S)=D(T)$. Then $S(0)=T(0)$, and the result of Theorem 3.5 holds.
\medskip

\noindent{\bf Proof.} Since $S$ and $T$ are self-adjoint with $D(S)=D(T)$,
it follows from [18, (ii) of Lemma 5.8] that $S(0)=T(0)$.
Hence, the result of Theorem 3.5 holds by Corollary 3.1. This completes the proof.
\medskip

To the end of this section, we shall consider the case that the perturbed relation $T$
can be  written as the following form:
$$T=S+A,                                                                                                                              \eqno(3.22)$$
where  $T, S, A\in LR(X)$ satisfy that
$$D(S)=D(T)=:D\subset D(A),                                                                                                       \eqno(3.23)$$
where  $S$ is the unperturbed relation and $A$ is the perturbed term.

We shall remark that in the single-valued case, any two operators $S$ and $T$  from $X$ into itself
with $D(S)=D(T)$ can be written as (3.22) with $A=T-S$. However, In the multi-valued case,
(3.22) may not hold since $S+(T-S)=T$ may not hold in general. It follows from Lemma 2.1
that $S+(T-S)=T$ holds, and then (3.22) holds with $A=T-S$,
 if and only if $D(T)\subset D(S)$ and $S(0)\subset T(0)$.

 In the following, we shall study what conditions $A$ satisfies such that $T$ is
  a trace class perturbation of  $S$.

In the case that $T, S, A\in CR(X)$, we established some relationships among their operator parts
$T_s$, $S_s$, and $A_s$ in [13]. The following result comes from  [13, Theorem 3.1].\medskip

\noindent {\bf Lemma 3.2.} Let  $T, S, A\in CR(X)$ satisfy (3.22) and (3.23).  Then
$$T_s=P_{T(0)^\perp}S_s+P_{T(0)^\perp}A_s\;\;{\rm in}\; D,                                                              \eqno(3.24)$$
where $P_{T(0)^\perp}$ is the orthogonal projection from $X$
onto $T(0)^\perp$. Furthermore, if $S$ and $A$ are Hermitian relations in $X^2$,
then $P_{T(0)^\perp}S_s$ and $P_{T(0)^\perp}A_s$
are Hermitian operators defined on $D$, respectively.\medskip

 Let  $T, S, A\in CR(X)$ satisfy (3.22) and (3.23). It follows from (2.3) and (3.23)  that
$$D(T_s)=D(S_s)=D,\;\; D(A_s)=D(A),                                                                                              \eqno(3.25)$$
$$R(T_s)\subset T(0)^\perp,\; R(S_s)\subset S(0)^\perp,\; R(A_s)\subset A(0)^\perp.                             \eqno(3.26)$$
 If $A(0)\subset S(0)$, then $S(0)=T(0)$ by (3.22), and consequently
 $P_{T(0)^\perp}S_s=S_s$ in $D$ by (3.26).
 Therefore, the following result directly follows from Lemma 3.2.
 \medskip

 \noindent {\bf Proposition 3.3.} Let  $T, S, A\in CR(X)$ satisfy  (3.22) and (3.23).
 If $A(0)\subset S(0)$, then
$$T_s=S_s+P_{S(0)^\perp}A_s\;\;{\rm in}\; D.                                                                            \eqno(3.27)$$

 \noindent {\bf Remark 3.2.} Proposition 3.3 is a generalization of [13, Corollary 3.1],
 where it is required that $A$ is single-valued.\medskip

\noindent{\bf Lemma 3.3.} Let  $S, T\in CR(X)$ satisfy that $S(0)=T(0)$, and $D(S)=D(T)=:D\subset S(0)^\perp$.
Then for every $\ld \in {\mathbf C}$,
$$(T-\lambda I)^{-1}(T-S)=(T_s-\lambda I)^{-1}(T_s-S_s).                                                        \eqno(3.28)  $$

\noindent{\bf Proof.} Fix any  $\ld \in {\mathbf C}$, and let
$$U_1=(T-\lambda I)^{-1}(T-S),\;\;U_2=(T_s-\lambda I)^{-1}(T_s-S_s).$$
Since $S(0)=T(0)$, by (2.1) and (2.3) we have that
$$T=T_s\oplus S_\infty, \;\;S=S_s\oplus S_\infty, \;T_\infty=S_\infty=\{0\}\times S(0),                                       \eqno (3.29)$$
$$D(T_s)=D(S_s)=D,\;R(T_s)\subset S(0)^\perp,\; R(S_s)\subset S(0)^\perp.                                      \eqno (3.30)$$

It follows from (3.29) that
$$U_1=(T-\lambda I)^{-1}\left((T_s-S_s)\oplus S_\infty\right).                                                                                        $$
So, for any given $(x,y)\in U_1$, there exists $z\in X$ such that
$$(x,z)\in (T_s-S_s)\oplus S_\infty,\;\; (z,y)\in (T-\lambda I)^{-1},                                                                          \eqno (3.31) $$
and then there exist $z_T\in R(T_s), z_S\in R(S_s),$ and $\omega \in S(0)$ such that
$$(x,z_T)\in T_s,\;\;(x,z_S)\in S_s,\;\; z=z_T-z_S+\omega.                                                                    \eqno (3.32) $$
It follows from the second relation in (3.31) that
 $(y,z+\ld y)\in T$, which, together with the first relation in (3.29)
and the third relation in (3.32), implies that
$$(y, z_T-z_S+\ld y+\omega)\in T_s\oplus S_\infty.                                                                           \eqno (3.33)  $$
Noting that $z_T\in R(T_s)\subset S(0)^\perp, z_S\in R(S_s)\subset S(0)^\perp$ by the last two relations
in (3.30), and $y\in D\subset S(0)^\perp$ by the assumption, we have that
$z_T-z_S+\ld y\in S(0)^\perp. $ Hence, it follows from (3.33) that
$$(y, z_T-z_S+\ld y)\in T_s,\;\; (0,\omega)\in  S_\infty.$$
The first relation in the above yields that $(y, z_T-z_S)\in T_s-\ld I$,
and so  $(z_T-z_S,y)\in (T_s-\ld I)^{-1}.$
In addition, it follows from the first two relations in (3.32) that $(x, z_T-z_S)\in T_s-S_s$.
Therefore, $(x,y)\in U_2$, and consequently $U_1\subset U_2$.

It can be easily verified that $U_2\subset U_1$  by the fact that $T_s\subset T$ and $S_s\subset S$.
Hence, $U_1=U_2$, and then (3.28) holds. This completes the proof.\medskip

\noindent {\bf Proposition 3.4.} Let  $S, T\in CR(X)$ satisfy that $S(0)=T(0)$, and $D(S)=D(T)=:D\subset S(0)^\perp$.
Then for every $\ld \in {\mathbf C}$,
$$(T-\lambda I)^{-1}-(S-\lambda I)^{-1}=-(T_s-\lambda I)^{-1}(T_s-S_s)(S-\lambda I)^{-1}.                    \eqno(3.34)  $$

\noindent{\bf Proof.} Fix any  $\ld \in {\mathbf C}$, and let $U:=(T-\lambda I)^{-1}-(S-\lambda I)^{-1}.$
By Lemma 2.2 we have that
$$U=-(T-\lambda I)^{-1}\left[(T-\ld I)-(S-\ld I)\right] (S-\lambda I)^{-1}.$$
It can be easily verified that $(T-\ld I)-(S-\ld I)=T-S$. So we get that
$$U=-(T-\lambda I)^{-1}(T-S) (S-\lambda I)^{-1},                                                                               $$
which implies that (3.34) holds by Lemma 3.3. The proof is complete.
\medskip

\noindent{\bf Theorem 3.6.} Let $T, S, A\in CR(X)$ be Hermitian and satisfy that
(3.22),  (3.23), $A(0)\subset S(0)$,  $S(0)^\perp \subset D(A)$, and
$\rho(S)\cap \rho(T)\ne \emptyset$.
If one of the following conditions is satisfied:
\begin{itemize}
\item[{\rm (i)}]  $P_{S(0)^\perp} A_s|_{S(0)^\perp}\in EB_1(S(0)^\perp)$;

\item[{\rm (ii)}]  $P_{S(0)^\perp} A_s|_{S(0)^\perp}$ is a finite rank operator on $S(0)^\perp$;

\end{itemize}
then $T$ is a trace class perturbation of $S$ in $X^2$.
\medskip

\noindent{\bf Proof.}  (i) Suppose that $P_{S(0)^\perp} A_s|_{S(0)^\perp}\in EB_1(S(0)^\perp)$.
By the assumption that $A(0)\subset S(0)$ one has that
$S(0)=T(0)$, and (3.27) holds by Proposition 3.3.
It follows from Lemma 2.4 that
$$ D\subset S(0)^\perp,\; R(T_s)\subset T(0)^\perp = S(0)^\perp,\;
\rho(T_s)=\rho(T),\;\rho(S_s)=\rho(S).                                                                             \eqno(3.35) $$
 Hence, $\rho(S_s)\cap \rho(T_s)=\rho(S)\cap \rho(T)\ne \emptyset$.
By Proposition 3.4 and (3.27) we have that for any $\ld \in \rho(S)\cap \rho(T)$,
$$\begin{array}{rrll}
&&(T-\lambda I)^{-1}-(S-\lambda I)^{-1}\\
&=&-(T_s-\lambda I)^{-1}(T_s-S_s)(S-\lambda I)^{-1}\\
&=&-(T_s-\lambda I)^{-1}P_{S(0)^\perp}A_s (S-\lambda I)^{-1}.
\end{array}                                                                                                                                       \eqno(3.36)  $$
Note that  $(S-\lambda I)^{-1}\in EB(X, S(0)^\perp)$, and
$(T_s-\lambda I)^{-1}\in EB(S(0)^\perp)$ by the first two relations in (3.35).
Since $P_{S(0)^\perp} A_s|_{S(0)^\perp}\in EB_1(S(0)^\perp)$,
we obtain that $(T-\lambda I)^{-1}-(S-\lambda I)^{-1}\in EB_1(X)$
by (3.36) and Lemma 2.6. Therefore, $T$ is a trace class perturbation of $S$ in $X^2$
by Theorem 3.4.

(ii)  Suppose that  $P_{S(0)^\perp} A_s|_{S(0)^\perp}$ is a finite rank operator on $S(0)^\perp$.
Then $P_{S(0)^\perp} A_s|_{S(0)^\perp}\in EB_1(S(0)^\perp)$ by their definitions.
Consequently,  $T$ is a trace class perturbation of $S$ in $X^2$ by the above assertion.
This completes the proof.\medskip

\noindent{\bf Lemma 3.4 {\rm [13, Theorem 3.2]}.} Let $S\in CR(X)$ be self-adjoint,
 $T, A\in CR(X)$ be Hermitian, and they satisfy  (3.22) and (3.23).
 Then (3.27) holds and
$$A(0)\subset S(0)=T(0).                                                                                                         \eqno(3.37)        $$

\noindent{\bf Corollary 3.3.} Let $S\in LR(X)$ be self-adjoint,  $T,A\in CR(X)$ be Hermitian,
and they satisfy that (3.22),  (3.23), and $S(0)^\perp \subset D(A)$.
If one of the conditions (i) and (ii) in Theorem 3.6 holds,
then $A(0)\subset S(0)$,  $T$ is self-adjoint in $X^2$,
and $T$ is a trace class perturbation of $S$ in $X^2$.
\medskip

\noindent{\bf Proof.} It follows from Lemmas 2.4 and 3.4 that
$D\subset S(0)^\perp$, $R(T_s)\subset S(0)^\perp$, and (3.27) and (3.37) hold.
By the assumption that one of the conditions (i) and (ii) in Theorem 3.6 holds,
 $P_{S(0)^\perp} A_s|_{S(0)^\perp}$ is compact, and then bounded on $S(0)^\perp$.
 Thus,  $P_{S(0)^\perp} A_s|_{S(0)^\perp}$ is $S_s$-bounded
 with $S_s$-bound $0$ by [20, Proposition  on Page 93].
So $T$ is self-adjoint in $X^2$ by [13, Theorem 4.2].
This implies that $\rho(S)\cap \rho(T)\ne \emptyset$. Hence,
all the assumptions of Theorem 3.6 hold, and then
$T$ is a trace class perturbation of $S$ in $X^2$ by Theorem 3.6.
The proof is complete.\medskip

The following result is a  direct consequence of Theorem 3.6.\medskip

\noindent{\bf Corollary 3.4.} Let $T, S, A\in CR(X)$ be Hermitian and satisfy that
(3.22),  (3.23),  $A(0)\subset S(0)$,   and $\rho(S)\cap \rho(T)\ne \emptyset$.
If one of the following conditions is satisfied:
\begin{itemize}
\item[{\rm (i)}]  $S(0)^\perp \subset D(A)$, and $ A_s|_{S(0)^\perp}\in EB_1(S(0)^\perp)$;

\item[{\rm (ii)}]  $S(0)^\perp \subset D(A)$, and $ A_s|_{S(0)^\perp}$ is a finite rank operator on $S(0)^\perp$;

\item[{\rm (iii)}] $ D(A)=X$, and $ A_s\in EB_1(X)$;

\item[{\rm (iv)}]  $ D(A)=X$, and $ A_s$ is a finite rank operator on $X$;

\end{itemize}
then $T$ is a trace class perturbation of $S$ in $X^2$.
\medskip

The following result directly follows from Corollary 3.3.\medskip

\noindent{\bf Corollary 3.5.} Let $S\in LR(X)$ be self-adjoint,  $T,A\in CR(X)$ be Hermitian,
and they satisfy that (3.22) and (3.23).
If one of the conditions (i) - (iv) in Corollary 3.4 is satisfied, then the results of Corollary 3.3 hold.\medskip

\noindent{\bf Remark 3.3.} The characterizations given in this section are very important in the study of 
problems about trace class perturbation of closed relations.
We shall apply them to study stability of absolutely continuous spectra of closed relations in Hilbert spaces under
trace class perturbation, and then discuss their applications to symmetric difference equations in our
forthcoming papers.

\bigskip

 \noindent{\bf \large References}
\def\hang{\hangindent\parindent}
\def\textindent#1{\indent\llap{#1\enspace}\ignorespaces}
\def\re{\par\hang\textindent}
\noindent \vskip 3mm
\re{[1]} R. Arens, Operational calculus of linear relations,
Pac. J. Math. 11(1961) 9--23.

\re{[2]} T. Ya. Azizov, J. Behrndt, P. Jonas, C. Trunk,
Compact and finite rank perturbations of linear
relations in Hilbert spaces, Integr. Equ. Oper. Theory 63 (2009) 151--163.

\re{[3]}  E. A. Coddington, Extension theory of formally normal
and symmetric subspaces, Mem. Am. Math. Soc. 134 (1973).

\re{[4]} R. Cross, Multivalued Linear Operators. In : Monographs and Textbooks in Pure and Applied Mathematics,
vol. 213. New York: Marcel Dekker; 1998.

\re{[5]} A. Dijksma, H. S. V. de Snoo, Eigenfunction extensions associated with pairs of ordinary differential expressions,
J. Differ. Equations 60(1985) 21--56.

\re{[6]} S. Hassi, H. de Snoo, One-dimensional graph perturbations of self-adjoint relations, Ann. Aca. Sci. Fenn.
Math. 20(1997) 123--164.

\re{[7]} S. Hassi, H. de Snoo, F. H. Szafraniec, Componentwise and cartesian decompositions of linear relations,
Dissertationes Math. 465, 2009 (59 pages).

\re{[8]} T. Kato, Perturbation Theory for Linear Operators, 2nd ed.,
Springer-Verlag, Berlin \newline/Heidelberg/New York/Tokyo,  1984.

\re{[9]} M. Lesch, M. Malamud, On the deficiency indices and self-adjointness of
symmetric Hamiltonian systems, J. Differ. Equations 18(2003) 556--615.

\re{[10]} J. Von Neumann, Functional Operator II: The Geometry of Orthogonal Spaces. Ann. Math. Stud. 22, Princeton U.P., 1950.

\re{[11]} M. Reed, B. Simon, Methods of Modern Mathematical Physics
I: Functional Analysis, Academic Press, 1972.

\re{[12]} G. Ren, Y. Shi, Defect indices and definiteness conditions for discrete linear Hamiltonian systems, Appl. Math. Comput.
218(2011) 3414--3429.

\re{[13]} Y. Shi, Stability of essential spectra of self-adjoint subspaces under compact perturbations,
 J. Math. Anal. Appl.  433(2016) 832--851.

\re{[14]} Y. Shi, C. Shao, G. Ren, Spectral properties of self-adjoint subspaces,
Linear Algebra Appl. 438(2013) 191--218.

\re{[15]} Y. Shi, H. Sun, Self-adjoint extensions for second-order
symmetric linear difference equations, Linear Algebra Appl. 434(2011) 903--930.

\re{[16]} Y. Shi, G. Xu, G. Ren, Boundedness and closedness of linear relations, Linear and Multilinear algebra 66(2018) 309--333.

\re{[17]} C. Tretter, Spectral Theory of Block Operator Matrices and Applications,
Imperial College Press, 2008.

\re{[18]} G. Xu, Y. Shi, Perturbations of spectra of closed subspaces in Banach spaces. Linear Algebra Appl.
     531(2017) 547--574.

\re{[19]} G. Xu, Y. Shi, Perturbations of essential spectra of self-adjoint relations under relatively
        compact perturbations, Linear and Multilinear algebra 66(12) (2018)  2438--2467.

\re{[20]} J. Weidmann, Linear Operators in Hilbert Spaces,
Graduate Texts in Math., vol.68, Springer-Verlag, New York/Berlin/Heidelberg/Tokyo, 1980.

\re{[21]} D. Wilcox, Essential spectral of linear relations, Linear Algebra Appl. 462(2014) 110--125.

\end{document}